\documentclass[12pt,a4paper,reqno]{amsart} 
\usepackage{amsmath,amssymb,amsthm}
\usepackage{enumerate} 
\usepackage{url} 
\usepackage{fullpage}  
  
\newcommand{\eps}{\varepsilon} 
\newcommand{\D}{\mathfrak{D}}
\newcommand{\M}{\mathfrak{M}}
\newcommand{\m}{\mathfrak{m}}
\newcommand{\tr}{\mathfrak{t}} 
\newcommand{\X}{\mathfrak{X}}

\newcommand{\dx}{\mathrm{d}}

\newcommand{\R}{\mathbb{R}}

\newcommand{\Odip}[2]{\mathcal{O}_{#1}\!\left(#2\right)\mathchoice{\!}{}{}{}}
\newcommand{\Odi}[1]{\Odip{}{#1}}
\newcommand{\odip}[2]{{o}_{#1}\!\left(#2\right)\mathchoice{\!}{}{}{}}
\newcommand{\odi}[1]{\odip{}{#1}}

\allowdisplaybreaks
 
\renewcommand{\qedsymbol}{$\square$}
\newenvironment{Proof}[1][Proof]{\par\noindent\textbf{#1.}~}
{\hfill\qedsymbol\smallskip\par}
\newtheorem{Theorem}{Theorem}
\newtheorem{Lemma}{Lemma} 

\newtheoremstyle{Nonumtheorems}
{10pt}
{6pt}
{\itshape}
{}
{\bfseries}
{.}
{.5em}
{\thmname{#1}\thmnote{ (#3)}}

\theoremstyle{Nonumtheorems}

\begin{document}

\title{On a ternary  Diophantine problem with mixed \\powers of primes} 
\author{Alessandro Languasco \& Alessandro Zaccagnini}
\date{}

 
\begin{abstract}
Let $1 < k < 33 / 29$.
We prove that if $\lambda_1$, $\lambda_2$ and $\lambda_3$ are non-zero
real numbers, not all of the same sign and that
$\lambda_1 / \lambda_2$ is irrational and $\varpi$ is any real number,
then for any $\eps > 0$ the inequality
$
  \bigl\vert
    \lambda_1 p_1 + \lambda_2 p_2^2 + \lambda_3 p_3^k
    +
    \varpi
  \bigr\vert
  \le
  \bigl( \max_j p_j \bigr)^{-(33 - 29 k) / (72 k) + \eps}
$
has infinitely many solutions in prime variables $p_1$, \dots, $p_k$ .  
\end{abstract}

\maketitle

\section{Introduction}
The goal of this paper is to solve a ternary Diophantine approximation
problem which involves real powers of prime numbers.
We restrict our attention to the values of the form
$ \lambda_1 p_1 + \lambda_2 p_2^2 + \lambda_3 p_3^k$, where $k\in \R$
and $k>1$, but similar cases can be attacked with this approach as well, 
see, e.g., \cite{LanguascoZ2012e}.
Our main result is 

\begin{Theorem}
\label{Th:appl2}
Let $1 < k < 33 / 29$ be a real number and assume that $\lambda_1$,
$\lambda_2$ and $\lambda_3$ are non-zero real numbers, not all of the
same sign and that $\lambda_1 / \lambda_2$ is irrational.
Let $\varpi$ be any real number. 
For any $\eps > 0$ the inequality
\begin{equation}
\label{main-ineq-bis}
  \bigl\vert 
    \lambda_1 p_1 + \lambda_2 p_2^2 + \lambda_3 p_3^k
    +
    \varpi
  \bigr\vert
  \le
  \bigl( \max_j p_j \bigr)^{-(33 - 29 k) / (72 k) + \eps}
\end{equation}
has infinitely many solutions in prime variables $p_1$, \dots, $p_k$ .
\end{Theorem}
The proof of Theorem \ref{Th:appl2} uses the variant of the circle method 
introduced by Davenport and Heilbronn to deal with Diophantine problems. 
Classical papers on this topic with integral $k$ are Vaughan's ones, 
\cite{Vaughan1974a,Vaughan1974b}, Baker and Harman \cite{BakerH1982}, 
Harman \cite{Harman1991}.
For non-integral $k$ we recall that Tolev \cite{Tolev1992} 
studied the values of the form $p_1^k + p_2^k + p_3^k$
and proved that, for every $k \in (1, 15 / 14)$, all sufficiently large 
real numbers $\varpi$ can be well approximated.

In order to deal with a problem with mixed non-integral powers, like
the present one, a key tool is a suitable estimate for the
$L^2$-norms of exponential sums over prime powers, see Theorems
\ref{k-root-BCP-Gallagher} and \ref{k-root-Saffari-Vaughan} of
\S\ref{mean-square-section}, which has some independent interest.
These results allow us to have a comparatively wide ``major arc''
while keeping the resulting error term under control.
This idea appeared in Br\"udern, Cook and Perelli \cite{BrudernCP1997}
and we exploit it also in \cite{LanguascoZ2012c} and
\cite{LanguascoZ2012e}.

\smallskip
\textbf{Acknowledgement.}
We thank Alberto Perelli for a discussion.

\section{Proof of Theorem \ref{Th:appl2}} 

We use the variant of the circle method introduced by Davenport and
Heilbronn to deal with Diophantine problems.
In order to prove that \eqref{main-ineq-bis} has infinitely many
solutions, it is sufficient to construct an increasing sequence $X_n$
with limit $+\infty$ such that \eqref{main-ineq-bis} has at least a
solution with $\max_j p_j \in [\delta X_n, X_n]$, where $\delta$ is a
small, fixed positive constant.
This sequence actually depends on rational approximations for
$\lambda_1 / \lambda_2$: more precisely, there are infinitely many
pairs of integers $a$ and $q$ such that $(a, q) = 1$, $q > 0$ and
\[
  \Bigl\vert \frac{\lambda_1}{\lambda_2} - \frac aq \Bigr\vert
  \le
  \frac 1{q^2}.
\]
We take the sequence $X = q^{9 k / (2 k + 3)}$ (dropping the useless
suffix $n$) and then, as customary, define all of the circle-method
parameters in terms of $X$.
We may obviously assume that $q$ is sufficiently large.
The choice of the exponent of $q$ and of all the other parameters is
justified in the discussion in \S\ref{sec:param}.

Let
\begin{equation}
\label{S-def}
S_{k}(\alpha)
 =
\sum_{X \leq p^{k} \leq  2X} \log p \ e(p^{k} \alpha ).
\end{equation}
As usual, we approximate to $S_k$ using the function
\[
  T_k(\alpha)
  =
  \int_{(\delta X)^{1/k}}^{X^{1/k}} e(t^k \alpha) \, \dx t
\]
and notice the simple inequality
\begin{equation}
\label{bd-Tkappa}
  T_k(\alpha)
  \ll_{\delta, k}
  X^{1 / k -1}
  \min\bigl( X, \vert \alpha \vert^{-1} \bigr).
\end{equation}
We detect solutions of \eqref{main-ineq-bis} by means of the function
$\widehat{K}_{\eta}(\alpha) = \max(0, \eta - \vert \alpha \vert)$
for $\eta > 0$, which, as the notation suggests, is the Fourier
transform of
\[
  K_{\eta}(\alpha)
  =
  \Bigl( \frac{\sin(\pi \eta \alpha)}{\pi \alpha} \Bigr)^2
\]
for $\alpha \ne 0$, and, by continuity, $K_{\eta}(0) = \eta^2$.
This relation transforms the problem of counting solutions of the
inequality \eqref{main-ineq-bis} into estimating suitable integrals.
We recall the trivial property
\begin{equation}
\label{bd-K(eta)-bis}
  K_{\eta}(\alpha)
  \ll
  \min ( \eta^2, \vert \alpha \vert^{-2} ).
\end{equation}

For any measurable subset $\X$ of $\R$ let
\[
  I(\eta, \varpi, \X)
  =
  \int_{\X} S_1(\lambda_1 \alpha) S_2(\lambda_2 \alpha)
    S_k(\lambda_3 \alpha) K_{\eta}(\alpha)
    e(\varpi \alpha) \, \dx \alpha.
\]
In practice, we take as $\X$ either an interval or a half line, or the
union of two such sets.
The starting point of the method is the observation that
\begin{align*}
  I(\eta, \varpi, \R)
  &=
  \sum_{p_1, p_2^2, p_3^k \in [\delta X, X]}
    \log p_1 \log p_2 \log p_3  
      \int_{\R}
    K_{\eta}(\alpha)
    e \bigl(
        (\lambda_1 p_1 + \lambda_2 p_2^2 + \lambda_3 p_3^k
        + \varpi) \alpha
      \bigr) \, \dx \alpha \\
  &=
  \sum_{p_1, p_2^2, p_3^k \in [\delta X, X]}
    \log p_1 \log p_2 \log p_3  
  \max(0,
       \eta - \vert \lambda_1 p_1 + \lambda_2 p_2^2 + \lambda_3 p_3^k + \varpi
              \vert) \\
  &\le
  \eta (\log X)^3 \mathcal{N}(X),
\end{align*}
where $\mathcal{N}(X)$ denotes the number of solutions of the
inequality \eqref{main-ineq-bis} with $p_1$, $p_2^2$,
$p_3^k \in [\delta X, X]$.
We now give the definitions that we need to set up the method.
More definitions will be given at appropriate places later.
We let $P = P(X) = X^{4 / (5 k) - \eps} $,
$\eta = \eta(X) = X^{- (33 - 29 k) / (72 k) + \eps}$, and
$R = R(X) = \eta^{-2} X^{(k - 1) / (4 k)} (\log X)^3$.
The choice for $P$ is justified at the end of \S\ref{subs:J4-bis}, the
one for $\eta$ at the end of \S\ref{sec:intermediate-bis} and the one
for $R$ at the end of \S\ref{sec:trivial-bis}.
See also the discussion in \S\ref{sec:param} for a more detailed
argument.
We now decompose $\R$ as $\M \cup \m \cup \tr$ where
\[
  \M
  =
  \Bigl[ -\frac PX, \frac PX \Bigr],
  \qquad
  \m
  =
  \Bigl( -R, -\frac PX \Bigr) \cup \Bigl(\frac PX, R \Bigr),
  \qquad
  \tr
  =
  \R \setminus( \M \cup \m),
\]
so that
\[
  I(\eta, \varpi, \R)
  =
  I(\eta, \varpi, \M)
  +
  I(\eta, \varpi, \m)
  +
  I(\eta, \varpi, \tr).
\]
These sets are called the major arc, the intermediate (or minor) arc
and the trivial arc respectively.
In \S\ref{sec:major-bis} we prove that the major arc yields the main term
for $I(\eta, \varpi, \R)$.
In order to show that the contribution of the intermediate arc does
not cancel the main term, we exploit the hypothesis that
$\lambda_1 / \lambda_2$ is irrational to prove that
$\vert S_1(\lambda_1 \alpha) \vert^{1 / 2}$ and
$\vert S_2(\lambda_2 \alpha) \vert$ can not both be large for
$\alpha \in \m$: see \S\ref{sec:intermediate-bis}, and in particular
Lemma~\ref{Lemma-approx-bis}, for the details.
The trivial arc, treated in \S\ref{sec:trivial-bis}, only gives a rather
small contribution.

From now on, implicit constants may depend on the coefficients
$\lambda_j$, on $k$, $\delta$ and $\varpi$.

\subsection{The major arc}
\label{sec:major-bis}

We write
\begin{align*}
  I(\eta, \varpi, \M)
  &=
  \int_{\M}
    S_1(\lambda_1 \alpha) S_2(\lambda_2 \alpha) S_k(\lambda_3 \alpha)
    K_{\eta}(\alpha) e(\varpi \alpha) \, \dx \alpha \\
  &=
  \int_{\M}
    T_1(\lambda_1 \alpha) T_2(\lambda_2 \alpha) T_k(\lambda_3 \alpha)
    K_{\eta}(\alpha) e(\varpi \alpha) \, \dx \alpha \\
  &\qquad+
  \int_{\M}
    \bigl( S_1(\lambda_1 \alpha) - T_1(\lambda_1 \alpha) \bigr)
    T_2(\lambda_2 \alpha) T_k(\lambda_3 \alpha)
    K_{\eta}(\alpha) e(\varpi \alpha) \, \dx \alpha \\
  &\qquad+
  \int_{\M}
    S_1(\lambda_1 \alpha)
    \bigl(S_2(\lambda_2 \alpha) - T_2(\lambda_2 \alpha) \bigr)
    T_k(\lambda_3 \alpha)
    K_{\eta}(\alpha) e(\varpi \alpha) \, \dx \alpha \\
  &\qquad+
  \int_{\M}
    S_1(\lambda_1 \alpha) S_2(\lambda_2 \alpha)
    \bigl(S_k(\lambda_3 \alpha) - T_k(\lambda_3 \alpha) \bigr)
    K_{\eta}(\alpha) e(\varpi \alpha) \, \dx \alpha \\
  &=
  J_1 + J_2 + J_3 + J_4,
\end{align*}
say.
We will give a lower bound for $J_1$ and upper bounds for $J_2$,
\dots, $J_4$.
For brevity, since the computations for $J_3$ are similar
to, but simpler than, the corresponding ones for $J_2$ and $J_4$, we
will skip them.

\subsection{Lower bound for $J_1$}
The lower bound $J_1 \gg \eta^2 X^{1 / 2 + 1 / k}$ is proved in a
classical way.
We have
\begin{align*}
  J_1
  &=
  \int_{\M}
    T_1(\lambda_1 \alpha) T_2(\lambda_2 \alpha) T_k(\lambda_3 \alpha)
    K_{\eta}(\alpha) e(\varpi \alpha) \, \dx \alpha \\
  &
  =
  \int_{\R}
    T_1(\lambda_1 \alpha) T_2(\lambda_2 \alpha) T_k(\lambda_3 \alpha)
     K_{\eta}(\alpha) e(\varpi \alpha) \, \dx \alpha \\
  &\qquad+
  \Odi{
  \int_{P / X}^{+\infty}
    \vert
      T_1(\lambda_1 \alpha) T_2(\lambda_2 \alpha) T_k(\lambda_3 \alpha)
    \vert K_{\eta}(\alpha) \, \dx \alpha}.
\end{align*}
Using inequalities \eqref{bd-Tkappa} and \eqref{bd-K(eta)-bis}, we see
that the error term is
\[
  \ll
  \eta^2
  X^{1 / k - 3 / 2}
  \int_{P / X}^{+\infty}
    \frac{\dx \alpha}{\alpha^3}
  \ll
  \eta^2 X^{1 / 2 + 1 / k} P^{-2}
  =
  \odi{\eta^2 X^{1 / 2 + 1 / k}}.
\]
For brevity, we set
$\D = [\delta X, X] \times [(\delta X)^{1/2}, X^{1/2}] \times
      [(\delta X)^{1/k}, X^{1/k}]$.
We can rewrite the main term in the form
\begin{align*}
  &
  \idotsint_{\D}
    \int_{\R}
      e \bigl(
      (\lambda_1 t_1 + \lambda_2 t_2^2 + \lambda_3 t_3^k
       + \varpi) \alpha
        \bigr) \, K_{\eta}(\alpha) \, \dx \alpha
        \, \dx t_1 \, \dx t_2 \, \dx t_3  \\
  &=
  \idotsint_{\D}
    \max(0, \eta
             -
             \vert \lambda_1 t_1 + \lambda_2 t_2^2 + \lambda_3 t_3^k
                    + \varpi
             \vert)
    \, \dx t_1 \, \dx t_2 \, \dx t_3.
\end{align*}
We now proceed to show that the last integral is $\gg \eta^2 X^{1/2+1/k}$.
Apart from trivial changes of sign, there are essentially three cases:

\begin{enumerate}

\item
$\lambda_1 > 0$, $\lambda_2 < 0$,  $\lambda_3 < 0$;

\item
$\lambda_1 > 0$, $\lambda_2 > 0$,  $\lambda_3 < 0$;

\item
$\lambda_1 > 0$, $\lambda_2 < 0$,  $\lambda_3 > 0$.

\end{enumerate}

We briefly deal with the second case, the other ones being similar.
A suitable change of variables shows that
\begin{align*}
  J_1
  &\gg
  \idotsint_{\D'}
    \max(0, \eta
            -
            \vert \lambda_1 u_1 + \lambda_2 u_2 + \lambda_3 u_3 
            \vert)
    \, \frac{\dx u_1 \, \dx u_2 \, \dx u_3}{u_2^{1/2} u_3^{1-1/k}} \\
  &\gg
  X^{1 / k - 3 / 2}
  \idotsint_{\D'}
    \max(0, \eta
            -
            \vert \lambda_1 u_1 + \lambda_2 u_2 + \lambda_3 u_3 
            \vert)
    \, \dx u_1 \, \dx u_2 \, \dx u_3 ,
\end{align*}
where $\D' = [\delta X, (1 - \delta) X]^3$, for large $X$.
For $j = 1$, $2$, let
$a_j = \vert \lambda_3 \vert \delta / \vert \lambda_j \vert$,
$b_j = 2 a_j$ and $\mathfrak{I}_j = [a_j X, b_j X]$.
Notice that if $u_j \in \mathfrak{I}_j$ for $j = 1$, $2$, then
\[
  \lambda_1 u_1 + \lambda_2 u_2 
  \in
    \bigl[
    2 \vert \lambda_3 \vert \delta X, 4 \vert \lambda_3 \vert \delta X
  \bigr]
\]
so that, for every such choice of $(u_1, u_2)$, the interval
$[a, b]$ with endpoints
$\pm \eta / \vert \lambda_3 \vert +
(\lambda_1 u_1 + \lambda_2 u_2) / \vert \lambda_3 \vert$
is contained in $[\delta X, (1 - \delta) X]$.
In other words, for $u_3 \in [a, b]$ the values of
$\lambda_1 u_1 + \lambda_2 u_2 + \lambda_3 u_3 $
cover the whole interval $[-\eta, \eta]$.
Hence, for any $(u_1, u_2) \in \mathfrak{I}_1 \times \mathfrak{I}_2$,
we have
\[
  \int_{\delta X}^{(1 - \delta) X}
    \hskip-4pt
    \max(0, \eta
            -
            \vert \lambda_1 u_1 + \lambda_2 u_2 + \lambda_3 u_3 
            \vert)
    \, \dx u_3
  =
  \frac1{\vert \lambda_3 \vert}
  \int_{-\eta}^{\eta} \max(0, \eta - \vert u \vert) \, \dx u
  \gg
  \eta^2.
\]
Finally,
\[
  J_1
  \gg
  \eta^2
  X^{1 / k - 3 / 2}
  \iint_{\mathfrak{I}_1 \times \mathfrak{I}_2}
    \dx u_1 \, \dx u_2
  \gg
  \eta^2 X^{1 / 2 + 1 / k},
\]
which is the required lower bound.

\subsection{Bound for $J_2$}

We define another approximation of $S_{k}(\alpha)$, namely
\begin{equation}
\label{U-def}  
\quad
  U_{k}(\alpha)
  =
    \sum_{X \leq n^{k} \leq  2X} e(n^{k} \alpha ).
\end{equation}
 
 The Euler summation formula implies that
\begin{equation}
\label{bd-T-U-bis}
  T_k(\alpha)
  -
  U_k(\alpha)
  \ll
  1 + \vert \alpha \vert X.
\end{equation}
Using \eqref{bd-K(eta)-bis} we see that
\begin{align*}
  J_2
  &\ll
  \eta^2
  \int_{\M}
    \bigl\vert S_1(\lambda_1 \alpha) - T_1(\lambda_1 \alpha) \bigr\vert \,
    \vert T_2(\lambda_2 \alpha) \vert \,
    \vert T_k(\lambda_3 \alpha) \vert \, \dx \alpha \\
  &\le
  \eta^2
  \int_{\M}
    \bigl\vert S_1(\lambda_1 \alpha) - U_1(\lambda_1 \alpha) \bigr\vert \,
    \vert T_2(\lambda_2 \alpha) \vert \,
    \vert T_k(\lambda_3 \alpha) \vert \, \dx \alpha \\
  &\qquad+
  \eta^2
  \int_{\M}
    \bigl\vert U_1(\lambda_1 \alpha) - T_1(\lambda_1 \alpha) \bigr\vert \,
    \vert T_2(\lambda_2 \alpha) \vert \,
    \vert T_k(\lambda_3 \alpha) \vert \, \dx \alpha \\
  &=
  \eta^2
  (A_2 + B_2),
\end{align*}
say.
In order to estimate $A_2$ we use Theorems \ref{k-root-BCP-Gallagher}
and \ref{k-root-Saffari-Vaughan}.
By the Cauchy inequality and \eqref{bd-Tkappa} above, for any fixed
$A > 0$ we have
\begin{align*}
  A_2
  &\ll
  \Bigl(
    \int_{-P / X}^{P / X}
      \bigl\vert S_1(\lambda_1 \alpha) - U_1(\lambda_1 \alpha) \bigr\vert^2
      \, \dx \alpha
  \Bigr)^{1 / 2}  
  \Bigl(
    \int_{-P / X}^{P / X}
      \vert T_2(\lambda_2 \alpha) \vert^2 \,
      \vert T_k(\lambda_3 \alpha) \vert^2 \,  \dx \alpha
  \Bigr)^{1 / 2} \\
  &\ll
    \Bigl(\frac X{(\log X)^A} \Bigr)^{1 / 2} 
    \Bigl(
    \int_0^{1 / X} X^{1 + 2 / k} \, \dx \alpha
    +
    \int_{1 / X}^{P / X} \frac{X^{2 / k - 3}}{\alpha^4} \dx \alpha
  \Bigr)^{1 / 2} \\
  &\ll_A
  \frac{X^{1 / 2 + 1 / k}}{(\log X)^{A / 2}}
\end{align*}
by Theorem \ref{k-root-Saffari-Vaughan} (with $C = 12 / 5$), which we
can use provided that $X / P \ge X^{1 / 6 + \eps} $, that is,
$P \le X^{5 / 6 - \eps}$.
This proves that $\eta^2 A_2 = \odi{\eta^2 X^{1 / 2 + 1 / k}}$.
Furthermore, using the inequalities \eqref{bd-Tkappa} and \eqref{bd-T-U-bis}
we see that
\begin{align*}
  B_2
  &\ll
  \int_0^{1 / X}
    \vert T_2(\lambda_2 \alpha) \vert \,
    \vert T_k(\lambda_3 \alpha) \vert \, \dx \alpha
  +
  X
  \int_{1 / X}^{P / X}
    \alpha \,
    \vert T_2(\lambda_2 \alpha) \vert \,
    \vert T_k(\lambda_3 \alpha) \vert \,  \dx \alpha \\
  &\ll
 X^{1 / k - 1 / 2}
  +
  X^{1 / k - 1 / 2}
  \int_{1 / X}^{P / X}  \, \frac{\dx \alpha}{\alpha} 
  \ll
  X^{1 / k - 1 / 2} \log P,
\end{align*}
so that $\eta^2 B_2 = \odi{\eta^2 X^{1 / 2 + 1 / k}}$.

\subsection{Bound for $J_4$}
\label{subs:J4-bis}

Inequality \eqref{bd-K(eta)-bis} implies that
\begin{align*}
  J_4
  &\ll
  \eta^2
  \int_{\M}
    \bigl\vert S_1(\lambda_1 \alpha) \bigr\vert \,
    \bigl\vert S_2(\lambda_2 \alpha) \bigr\vert \,
    \bigl\vert S_k(\lambda_3 \alpha) - T_k(\lambda_3 \alpha) \bigr\vert
    \, \dx \alpha \\
  &\ll
  \eta^2
  \int_{\M}
    \bigl\vert S_1(\lambda_1 \alpha) \bigr\vert \,
    \bigl\vert S_2(\lambda_2 \alpha) \bigr\vert \,
    \bigl\vert S_k(\lambda_3 \alpha) - U_k(\lambda_3 \alpha) \bigr\vert
    \, \dx \alpha \\
  &\qquad+
  \eta^2
  \int_{\M}
    \bigl\vert S_1(\lambda_1 \alpha) \bigr\vert \,
    \bigl\vert S_2(\lambda_2 \alpha) \bigr\vert \,
    \bigl\vert U_k(\lambda_3 \alpha) - T_k(\lambda_3 \alpha) \bigr\vert
    \, \dx \alpha \\
  &=
  \eta^2 (A_4 + B_4),
\end{align*}
say.
The Parseval inequality and trivial bounds yield, for any fixed
$A > 0$,
\begin{align*}
  A_4
  &\ll
  X^{1 / 2}
  \Bigl(
    \int_{\M} \bigl\vert S_1(\lambda_1 \alpha) \bigr\vert^2 \, \dx \alpha
  \Bigr)^{1 / 2}
  \Bigl(
    \int_{\M}
      \bigl\vert S_k(\lambda_3 \alpha) - U_k(\lambda_3 \alpha) \bigr\vert^2
      \, \dx \alpha
  \Bigr)^{1 / 2} \\
  &\ll
  X (\log X)^{1 / 2}
  \frac PX J_k\Bigl(X, \frac XP \Bigr)^{1 / 2}
  \ll_A
  X^{1 / 2 + 1 / k} (\log X)^{1 / 2 - A / 2}
\end{align*}
by Theorems \ref{k-root-BCP-Gallagher} and \ref{k-root-Saffari-Vaughan}
(with $C = 12 / 5$) which we can use provided that
$X / P \ge X^{1 - 5 / (6k) + \eps}$, that is,
$P \le X^{5 / (6k) - \eps}$.
This proves that $\eta^2 A_4 = \odi{\eta^2 X^{1 / 2 + 1 / k}}$.
Furthermore, using \eqref{bd-T-U-bis}, the H\"older inequality and
trivial bounds we see that
\begin{align*}
  B_4
  &\ll
  \int_0^{1 / X}
    \bigl\vert S_1(\lambda_1 \alpha) \bigr\vert \,
    \bigl\vert S_2(\lambda_2 \alpha) \bigr\vert \,  \dx \alpha 
  +
  X
  \int_{1 / X}^{P / X}
    \alpha
    \bigl\vert S_1(\lambda_1 \alpha) \bigr\vert \,
    \bigl\vert S_2(\lambda_2 \alpha) \bigr\vert \,  \dx \alpha \\
  &\ll
  X^{1 / 2}
  +
  X
  \Bigl(
    \int_{1 / X}^{P / X}
      \bigl\vert S_1(\lambda_1 \alpha) \bigr\vert^2 \, \dx \alpha
  \Bigr)^{1 / 2}  
  \Bigl(
    \int_{1 / X}^{P / X}
      \bigl\vert S_2(\lambda_2 \alpha) \bigr\vert^4 \, \dx \alpha
    \int_{1 / X}^{P / X}
      \alpha^4 \, \dx \alpha
  \Bigr)^{1 / 4} \\
  &\ll
  X (X \log X)^{1 / 2}
  (X (\log X)^2)^{1 / 4}
  \Bigl( \frac PX \Bigr)^{5 / 4}
  \ll
  P^{5 / 4} X^{1 / 2} \log X.
\end{align*}
Here we used Satz~3 of Rieger \cite{Rieger1968} to bound the fourth
moment of $S_2$.
Hence, taking $P = \odi{X^{4 / (5 k)} (\log X)^{-1}}$ we get
$\eta^2 B_4 = \odi{\eta^2 X^{1 / 2 + 1 / k}}$.
We may therefore choose
\begin{equation}
\label{P-choice}
  P
  =
  X^{4 / (5 k) - \eps}.
\end{equation}

\subsection{The intermediate arc}
\label{sec:intermediate-bis}

We need to show that $\vert S_1(\lambda_1 \alpha) \vert^{1 / 2}$ and
$\vert S_2(\lambda_2 \alpha) \vert$ can not both be large for
$\alpha \in \m$, exploiting the fact that $\lambda_1 / \lambda_2$ is
irrational.
We do this using two famous results by Vaughan and Ghosh,
respectively, about $S_{1}(\alpha)$ and $S_2(\alpha)$.

\begin{Lemma}[Vaughan \cite{Vaughan1997}, Theorem 3.1]
\label{Vaughan-estim-bis}
Let $\alpha$ be a real number and $a,q$ be positive integers
satisfying $(a, q) = 1$ and $\vert \alpha - a / q \vert < q^{-2}$.
Then
\[
  S_{1}(\alpha)
  \ll
  \Bigl(\frac{X}{\sqrt{q}} + \sqrt{Xq} + X^{4/5} \Bigr) \log^4 X.
\]
\end{Lemma}
 
\begin{Lemma}[Ghosh \cite{Ghosh1981}, Theorem 2]
\label{Ghosh-estim}
Let $\alpha$ be a real number and $a,q$ be positive integers
satisfying $(a, q) = 1$ and $\vert \alpha - a / q \vert < q^{-2}$.
Let moreover $\eps > 0$.
Then
\[
  S_{2}(\alpha)
  \ll_{\eps}
  X^{1/2 + \eps}
  \Bigl(
    \frac{1}{q}
    +
    \frac{1}{X^{1/4}}
    +
    \frac{q}{X}
  \Bigr)^{1/4}.
\]
\end{Lemma}

\begin{Lemma}
\label{Lemma-approx-bis}
Let $1 \le k < 33 / 29$.
Assume that $\lambda_1 / \lambda_2$ is irrational and let
$X = q^{9 k / (2 k + 3)}$, where $q$ is the denominator of a
convergent of the continued fraction for $\lambda_1 / \lambda_2$.
Let $V(\alpha) =
\min \bigl( \vert S_1(\lambda_1 \alpha) \vert^{1 / 2},
            \vert S_2(\lambda_2 \alpha) \vert \bigr)$.
Then we have
\[
  \sup_{\alpha \in \m} V(\alpha)
  \ll
  X^{(29 k + 3) / (72 k) + \eps}.
\]
\end{Lemma}

\begin{Proof}
Let $\alpha \in \m$ and $Q = X^{(7 k - 3) / (18 k)} \leq P$.
By Dirichlet's Theorem, there exist integers $a_{i},q_{i}$  with
$1\leq q_{i}\leq X/Q$ and $(a_{i},q_{i})=1$, such that
$\vert \lambda_{i} \alpha q_{i}-a_{i}\vert \leq Q/X$, for $i=1,2$.
We remark that $a_{1}a_{2} \neq 0$, for otherwise we would have
$\alpha \in \M$.
Now suppose that $q_{i} \leq Q$ for $i=1,2$.
In this case we get
\[
  a_{2}q_{1} \frac{\lambda_{1}}{\lambda_{2}} - a_{1}q_{2}
  =
  ( \lambda_{1} \alpha q_{1}-a_{1}) \frac{a_{2}}{\lambda_{2} \alpha}
  -
  ( \lambda_{2} \alpha q_{2}-a_{2}) \frac{a_{1}}{\lambda_{2} \alpha}
\]
and hence
\begin{equation}
\label{bd-1-bis}
  \left\vert
    a_{2}q_{1} \frac{\lambda_{1}}{\lambda_{2}} - a_{1}q_{2}
  \right\vert
  \leq
  2\left(
    1 + \left\vert  \frac{\lambda_{1}}{\lambda_{2}} \right\vert
  \right)
  \frac{Q^{2}}{X}
  <
  \frac{1}{2q}
\end{equation}
for sufficiently large $X$.
Then, from the law of best approximation and the definition of $\m$,
we obtain
\begin{equation}
\label{bd-2-bis}
  X^{(2 k + 3) / (9 k)}
  =
  q
  \leq
  \vert a_{2}q_{1} \vert
  \ll
  q_{1} q_{2} R
  \leq
  Q^2 R
  \leq
  X^{(2 k + 3) / (9 k) - \eps},
\end{equation}
which is absurd.
Hence either $q_{1}>Q$ or $q_{2}>Q$.
Assume that $q_{1}>Q$.
Using Lemma \ref{Vaughan-estim-bis} on $S_1(\lambda_1 \alpha)$, we have
\begin{align*}
 V(\alpha)
 \leq
 \vert S_1(\lambda_1 \alpha) \vert^{1 / 2}
 &\ll
 \sup_{Q < q_{1} \leq X / Q}
 \left(
   \frac{X}{\sqrt{q_{1}}}
   +
   \sqrt{Xq_{1}}
   +
   X^{4/5}
  \right)^{1 / 2}
  \log^{2}X \\
  &\ll
  X^{(29 k + 3) / (72 k)}
  (\log X)^{2}.
\end{align*}
The other case is similar, using Lemma~\ref{Ghosh-estim} instead, and
hence Lemma \ref{Lemma-approx-bis} follows.
\end{Proof}

\begin{Lemma}
\label{Lemma:bd-minor-bis}
For $j = 1$, $2$ we have
\begin{align*}
  \int_{\m}
    \vert S_j(\lambda_j \alpha) \vert^{2 j} K_{\eta}(\alpha) \, \dx \alpha
  &\ll
  \eta X (\log X)^j \\
  \int_{\m}
    \vert S_k(\lambda_3 \alpha) \vert^2 K_{\eta}(\alpha) \, \dx \alpha
  &\ll
  \eta X^{1/k} (\log X)^3.
\end{align*}
\end{Lemma}

\begin{Proof}
The proof is achieved arguing as in \S\ref{sec:trivial-bis} below
where we bound the quantities $A$, $B$ and $C$, the main difference
being the fact that we have to split the range $[P / X, R]$ into two
intervals in order to use \eqref{bd-K(eta)-bis} efficiently.
See also the proof of Lemma~7 of Tolev \cite{Tolev1992}.
For the sake of brevity we skip the details.
\end{Proof}

Now let
\begin{align*}
  \X_1
  &=
  \{ \alpha \in [P / X, R] \colon
    \vert S_1(\lambda_1 \alpha) \vert^{1 / 2}
    \le
    \vert S_2(\lambda_2 \alpha) \vert \} \\
  \X_2
  &=
  \{ \alpha \in [P / X, R] \colon
    \vert S_1(\lambda_1 \alpha) \vert^{1 / 2}
    \ge
    \vert S_2(\lambda_2 \alpha) \vert \}
\end{align*}
so that $[P / X, R] = \X_1 \cup \X_2$ and
\[
  \Bigl\vert I(\eta, \varpi, \m) \Bigr\vert
  \ll
  \Bigl( \int_{\X_1} + \int_{\X_2} \Bigr)
    \bigl\vert
      S_1(\lambda_1 \alpha) S_2(\lambda_2 \alpha)
      S_k(\lambda_3 \alpha)  
    \bigr\vert
    K_{\eta}(\alpha) \, \dx \alpha.
\]
H\"older's inequality gives
\begin{align*}
  \int_{\X_1}
  &\le
  \max_{\alpha \in \X_1} \vert S_1(\lambda_1 \alpha) \vert^{1 / 2}
  \Bigl(
    \int_{\X_1}
      \vert S_1(\lambda_1 \alpha) \vert^2 K_{\eta}(\alpha) \, \dx \alpha
  \Bigr)^{1 / 4} \times \\
  &\qquad
  \Bigl(
    \int_{\X_1}
      \vert S_2(\lambda_2 \alpha) \vert^4 K_{\eta}(\alpha) \, \dx \alpha
  \Bigr)^{1 / 4}  
  \Bigl(
    \int_{\X_1}
      \vert S_k(\lambda_3 \alpha) \vert^2 K_{\eta}(\alpha) \, \dx \alpha
  \Bigr)^{1 / 2} \\ 
  &\ll
  \eta X^{(65 k + 39) / (72 k) + \eps}
\end{align*}
by Lemmas~\ref{Lemma-approx-bis} and \ref{Lemma:bd-minor-bis}.
The computation on $\X_2$ is similar and gives the same final result.  
Summing up,
\[
  \Bigl\vert I(\eta, \varpi, \m) \Bigr\vert
  \ll
  \eta X^{(65 k + 39) / (72 k) + \eps}
\]
and this is $\odi{\eta^2 X^{1 /2 + 1 / k}}$ provided that
\begin{equation}
\label{eta-choice}
  \eta
  =
  \infty\bigl( X^{(29 k - 33) / (72 k) + \eps} \bigr).
\end{equation}

\subsection{The trivial arc}
\label{sec:trivial-bis}

Using the H\"older inequality and a trivial bound for
$S_k(\lambda_3 \alpha)$ we see that
\begin{align*}
  \Bigl\vert I(\eta, \varpi, \tr) \Bigr\vert
  &\le
  2
  \int_R^{+\infty}
    \vert S_1(\lambda_1 \alpha) \vert \, \vert S_2(\lambda_2 \alpha) \vert \,
    \vert S_k(\lambda_3 \alpha) \vert \, 
    K_{\eta}(\alpha) \, \dx \alpha \\
  &\ll
  \Bigl(
    \int_R^{+\infty}
      \vert S_1(\lambda_1 \alpha) \vert^2 \, K_{\eta}(\alpha) \, \dx \alpha
  \Bigr)^{1 / 2} 
  \Bigl(
    \int_R^{+\infty}
      \vert S_2(\lambda_2 \alpha) \vert^4 \, K_{\eta}(\alpha) \, \dx \alpha
  \Bigr)^{1 / 4} \times \\
  &\qquad
  \Bigl(
    \int_R^{+\infty}
      \vert S_k(\lambda_3 \alpha) \vert^4 \, K_{\eta}(\alpha) \, \dx \alpha
  \Bigr)^{1 / 4} \\
  &\le
  X^{1 / 2 k}
  \Bigl(
    \int_R^{+\infty}
      \vert S_1(\lambda_1 \alpha) \vert^2 \, K_{\eta}(\alpha) \, \dx \alpha
  \Bigr)^{1 / 2}  
  \Bigl(
    \int_R^{+\infty}
      \vert S_2(\lambda_2 \alpha) \vert^4 \, K_{\eta}(\alpha) \, \dx \alpha
  \Bigr)^{1 / 4} \times \\
  &\qquad
  \Bigl(
    \int_R^{+\infty}
      \vert S_k(\lambda_3 \alpha) \vert^2 \, K_{\eta}(\alpha) \, \dx \alpha
  \Bigr)^{1 / 4} \\
  &\ll
  X^{1 / 2 k}
  A^{1 / 2}
  B^{1 / 4}
  C^{1 / 4},
\end{align*}
say, where in the last but one line we used the inequality
\eqref{bd-K(eta)-bis}, and we set
\[
  A
  =
  \int_{\vert \lambda_1 \vert R}^{+\infty}
    \frac{\vert S_1(\alpha) \vert^2}{\alpha^2} \, \dx \alpha,
  \quad
  B
  =
  \int_{\vert \lambda_2 \vert R}^{+\infty}
    \frac{\vert S_2(\alpha) \vert^4}{\alpha^2} \, \dx \alpha,
  \quad
  C
  =
  \int_{\vert \lambda_3 \vert R}^{+\infty}
    \frac{\vert S_k(\alpha) \vert^2}{\alpha^2} \, \dx \alpha.
\]
We have
\[
  A
  \ll
  \sum_{n \ge \vert \lambda_1 \vert R}
    \frac 1{(n - 1)^2}
    \int_{n - 1}^n \vert S_1(\alpha) \vert^2 \, \dx \alpha
  \ll
  \frac{X \log X}{\vert \lambda_1 \vert R}
\]
by the Prime Number Theorem (PNT).
Arguing similarly, using again Satz~3 of Rieger \cite{Rieger1968} and
Lemma~7 of Tolev \cite{Tolev1992} respectively, we see that we also
have $B \ll X (\log X)^2 / R$ and $C \ll X^{1 / k} (\log X)^3 / R$.
Collecting these estimates, we conclude that
\[
  \Bigl\vert I(\eta, \varpi, \tr) \Bigr\vert
  \ll
  \frac{X^{3 / 4 + 3 / (4 k)} (\log X)^2} R.
\]
Hence, $\bigl\vert I(\eta, \varpi, \tr) \bigr\vert =
  \odi{\eta^2 X^{1 / 2 + 1 / k}}$
provided that we choose, say,
\begin{equation}
\label{R-choice}
  R
  =
  \eta^{-2} X^{(1 - 1 / k) / 4} (\log X)^3.
\end{equation}

\subsection{Remark on the choice of the parameters}
\label{sec:param}

The constraint on the choice $X = q^{9 k / (2 k + 3)}$ with
$1 < k < 33 / 29$ arises from the bounds \eqref{bd-1-bis} and
\eqref{bd-2-bis}.
Their combination prevents us from choosing the optimal value
$X = q^2$.
This is justified as follows: neglecting log-powers, let
$X = q^{a(k)}$, $Q = X^{b(k)}$, $\eta = X^{-c(k)}$, and recall the
choices $P = X^{4 / (5 k) - \eps}$ in \eqref{P-choice} and
$R = \eta^{-2} X^{(1 - 1 / k) / 4} (\log X)^3$ in \eqref{R-choice}
which are due, respectively, to the bound for $B_4$ and for the
trivial arc.
Then, essentially, we have to maximize $k$ subject to the constraints
\[
  \begin{cases}
    a(k) \ge 1 \\
    0 \le b(k) \le \frac 4{5 k} \\
    c(k) \ge 0 \\
    2 b(k) - 1 \le - 1 / a(k)
    &\text{by \eqref{bd-1-bis},} \\
    2 b(k) + 2 c(k) + \frac14 (1 - \frac1k) \le 1 / a(k)
    &\text{by \eqref{bd-2-bis},} \\
    - c(k) \ge \frac12  - \frac1{2 k} - \frac14 b(k)
    &\text{by \eqref{eta-choice},}
  \end{cases}
\]
which is a linear optimization problem in the variables $1 / a(k)$,
$b(k)$, $c(k)$ and $1 / k$.
The solution for this problem is $1 / a(k) = (2 k + 3) / (9 k)$,
$b(k) = (7 k - 3) / (18 k)$, $c(k) = (33 - 29 k) / (72 k)$,
for $1 / k \ge 29 / 33$, and
this is equivalent to the statement of Theorem~\ref{Th:appl2}.

\section{$L^2$-norms of exponential sums over prime powers}
\label{mean-square-section}

In the proof of Theorem \ref{Th:appl2} we needed a mean-square
average of $S_{k}(\alpha) - U_{k}(\alpha)$, respectively defined
in \eqref{S-def} and \eqref{U-def}, for $k> 1$.  
In this section we see the slightly more general case
$k > 0$.

We need to recall that 
$\theta(x)=\sum_{p\leq x} \log p$ and to define the quantity
\begin{equation}
\label{k-root-Selberg-int-def}
J_{k}(X,h) 
=
\int_{X}^{2X}
\Bigl(
\theta((x+h)^{1/k})- \theta(x^{1/k}) -((x+h)^{1/k}-x^{1/k})
\Bigr)^2
\dx x
\end{equation}
which is a generalization of the Selberg integral which is well suited for
 our problem.
To be consistent with the classical definition, we will also denote
$J_1$ as $J$.

We want first to relate a truncated $L^2$-average 
of $ S_{k}(\alpha) - U_{k}(\alpha)$  with  $J_{k}(X,h)$
and then to obtain a suitable estimate for the latter.
\begin{Theorem} 
\label{k-root-BCP-Gallagher}
Let $k > 0$ be a real number.
For $0 < Y \leq 1/2$ we have
\[
\int_{-Y}^Y
\vert
S_{k}(\alpha) - U_{k}(\alpha) 
\vert^2
\dx \alpha
\ll_{k} 
  \frac{X^{2/k-2}\log^{2}X}{Y} +
Y^2X
+
Y^2 J_{k} \Bigl( X,\frac{1}{2Y} \Bigr),
\]
where $J_{k}(X,h)$ is defined in \eqref{k-root-Selberg-int-def}.
\end{Theorem}

A similar result holds replacing $\log p$ with the von Mangoldt
function $\Lambda(n)$ in the definition of $S_{k}(\alpha)$ in
\eqref{S-def}; the only difference in the statement above will be
replacing $J_k$ with $J_{k,\psi}$ as defined in \eqref{def_Jpsi}.
The case  $k=1$ is well known, see, e.g., Lemma 1 of
Br\"udern-Cook-Perelli \cite{BrudernCP1997}.

In order to state the following result, we introduce an hypothesis on
the density of the zeros of the Riemann zeta-function.
With classical notation, we assume that there exist constants $B \ge 0$
and $C \ge 2$ such that for $\sigma \in \bigl[ 1/2, 1 \bigr]$ and
$T \ge 2$ we have
\begin{equation}
\label{I-H}
  N(\sigma, T)
  \ll
  T^{C (1 - \sigma)} (\log T)^{B}.
\end{equation}
Huxley \cite{Huxley1972} proved that \eqref{I-H} holds with $C = 12/5$
and some $B \ge 0$.
\begin{Theorem}
\label{k-root-Saffari-Vaughan}
Let $k > 0$ be a real number and $\eps$ be an arbitrarily small
positive constant.
Assuming that \eqref{I-H} holds,
there exists a positive constant $c_1 = c_{1}(\eps)$, which does not
depend on $k$, such that
\[
J_{k}(X,h)
\ll_{k}
h^2
X^{2/k-1}
\exp
\Big(
- c_{1} 
\Big(
\frac{\log X}{\log \log X}
\Big)^{1/3}
\Big)
\]
uniformly for  $X^{1-2/(Ck)+\eps} \leq h \leq X$.  
Assuming further that the Riemann Hypothesis holds, we have
\[
J_{k}(X,h)
\ll_{k} 
h X^{1/k}  \log^2 \Bigl(\frac{2X}{h}\Bigr)  
\]
uniformly for $ X^{1-1/k} \leq h  \leq X$.
\end{Theorem}

Notice that if  $k\geq 1$ and $h \le X^{1 - 1 / k}$, then the bound $J_k(X, h) \ll X \log X$
follows immediately from the Prime Number Theorem. For $k<1$ the previous
condition on $h$ become essentially meaningless.
The case $k=1$ of the previous theorem was proved by
Saffari-Vaughan \cite{SaffariV1977}, see \S6 there,
while
in Zaccagnini \cite{Zaccagnini1998} a wider range for $h$ in the
unconditional case is given, but the proof does not easily lend itself
to the generalization we pursue here.
The unconditional case $k = 2$ of Theorems \ref{k-root-BCP-Gallagher}
and \ref{k-root-Saffari-Vaughan} (with $C=12/5$) was proved in
Languasco-Settimi \cite{LanguascoS2012}.
Applications of this case to some diophantine problems with primes and
squares of primes were given in \cite{LanguascoS2012} and in
Languasco-Zaccagnini \cite{LanguascoZ2012c}.
Results similar to Theorem \ref{k-root-Saffari-Vaughan} hold also
replacing $\theta(x)$ with $\psi(x)$, see Lemma
\ref{Lemma_da_h_a_deltax} below, and replacing $h$ with $\delta x$,
see Lemma \ref{Lemma_deltax}.

\subsection{Proof of Theorem \ref{k-root-BCP-Gallagher}}

Letting
${\mathcal I} := \int_{-Y}^Y
\vert
S_{k}(\alpha) - U_{k}(\alpha) 
\vert^2
\, \dx \alpha$,   
we see that the result is trivial for $0 < Y < 1 / X$ since
${\mathcal I} \ll Y X^{2/k} \leq Y^{-1}X^{2/k-2}\log^{2}X$ in this range.
Assuming that $1 / X \le Y \le 1 / 2$, we can write
\begin{align*}
{\mathcal I} 
&= 
 \int_{-Y}^{Y}
\Big \vert
 \sum_{X \leq n^{k} \leq  2X} ( \ell(n)-1 )e(n^{k} \alpha)
\Big\vert^{2}
\dx \alpha, 
\end{align*}
where $\ell(n) = \log p$ if $n=p$ prime and $ \ell(n)=0$ otherwise.  
By Gallagher's lemma (Lemma 1 of \cite{Gallagher1970})
we obtain
\[
{\mathcal I}  
\ll  
Y^{2}
\int_{-\infty}^{\infty}
\Bigl(
\sum_{\substack{x \leq n^{k} \leq x + H \\ X \leq n^{k} \leq  2X}}
 ( \ell(n)-1)
\Bigr)^{2}
\dx x  
\]
where we defined 
$H=1/(2Y)$.
We can restrict the integration range
to 
$
E = \left[  X - H, 2X\right]
$ since otherwise the inner sum is empty.
Moreover we split $E$ as $E=E_{1} \cup E_{2} \cup E_{3}$ where 
$E_{1}= \left [X - H,   X \right]$,
$E_{2}= \left [ X, 2X - H\right]$,
$E_{3}= \left [2X - H, 2X\right]$.
Accordingly  we can write
\begin{equation}
\label{I123}
\begin{split}
{\mathcal I}  & \ll 
Y^{2}
\left(
\int_{E_{1}}
+
\int_{E_{2}}
+
\int_{E_{3}}
\right)
\Bigl(
\sum_{\substack{x \leq n^{k} \leq x + H \\ X \leq n^{k} \leq  2X}}
 ( \ell(n)-1)
\Bigr)^{2}
\dx x
=
Y^{2}
(I_{1}+I_{2}+I_{3}),
\end{split}
\end{equation}
say.
We now proceed to estimate $I_{i}$, for every $i=1,2,3$.

\paragraph{\textbf{Estimation of $I_1$.}}
We immediately have
\begin{align}
I_{1} 
& \ll 
\notag
\int_{X - H}^{X}
\Bigl(
\theta\left((x+H)^{1/k}\right) - \theta(X^{1/k}) - \left((x + H)^{1/k} - X^{1/k} \right)
\Bigr)^{2}
\dx x 
+ H.
\end{align}
By trivial estimates  we obtain
\begin{align}
\label{I1}
I_{1}   
\ll
\log^{2}X
\int_{X - H}^{X}
\Bigl(
(x+H)^{1/k} - X^{1/k}
\Bigr)^{2}
\ \dx x 
+ H 
\ll
H^{3}X^{2/k-2}\log^{2} X
+ H.
\end{align}

\paragraph{\textbf{Estimation of $I_3$.}}
A similar argument gives the same bound for $I_3$, too: we omit it for
brevity.

\paragraph{\textbf{Estimation of $I_2$.}} 
We have
\begin{align}
\label{I2}
I_{2} 
& \ll
\int_{X }^{2X}
\Bigl(
\theta \left( (x + H)^{1/k} \right) - \theta \left( x^{1/k} \right) 
-
\left( (x + H)^{1/k}  - x^{1/k}\right)
\Bigr)^{2}
\dx x
+ 
{X} 
\notag \\
& =
J_{k}\left(X,H\right) +X,
\end{align}
where we used the definition in \eqref{k-root-Selberg-int-def}. 
Therefore, by \eqref{I123}-\eqref{I2}, the bound $Y\geq 1/X$ 
and recalling that $H=1/(2Y)$, we have
\begin{align*}
{\mathcal I}
&
\ll 
  \frac{X^{2/k-2}\log^{2}X}{Y} +  XY^{2} + Y^{2} J_{k}\left(X,\frac{1}{2Y}\right)
\end{align*} 
and this proves Theorem \ref{k-root-BCP-Gallagher}.

\subsection{Proof of Theorem \ref{k-root-Saffari-Vaughan}}

We reduce our problem to estimate
\begin{equation}
\label{def_Jpsi}
J_{k,\psi}(X,h)
:=
\int_{X}^{2X}
\left(
\psi((x+h)^{1/k})-\psi(x^{1/k}) -
((x+h)^{1/k}-x^{1/k})
\right)^2
\dx x 
\end{equation}
since, using $|a+b|^{2} \leq 2|a|^{2}+2|b|^{2}$ and Lemma
\ref{psi-theta-passaggio} below, we have
\begin{align}
\notag 
J_{k}(X,h) 
&\ll
J_{k,\psi}(X,h)  
+
\int_{X}^{2X}
\left(
\psi((x+h)^{1/k})- \psi(x^{1/k}) -
\theta((x+h)^{1/k})+\theta(x^{1/k}) 
\right)^{2}
\dx x \\ 
\label{fromThetatoPsiFinale}
&
\ll
J_{k,\psi}(X,h) 
+
h X^{1/k}.
\end{align} 
To estimate the right-hand side of \eqref{fromThetatoPsiFinale}, we
use the following result we will prove later.

\begin{Lemma} 
\label{Lemma_da_h_a_deltax}
Let $k >0$ be a real number and 
$\eps$ be an arbitrarily small positive constant.  
Assuming that \eqref{I-H} holds,
there exists a positive constant $c_1 = c_{1}(\eps)$, which does not
depend on $k$,  such that
\[
J_{k,\psi}(X,h)
\ll
h^2
X^{2/k-1}
\exp
\Big(
- c_{1}
\Big(
\frac{\log X}{\log \log X}
\Big)^{1/3}
\Big)
\]
uniformly for  $X^{1-2/(Ck)+\eps} \leq h \leq X$, 
where $J_{k,\psi}(X,h)$ is defined in \eqref{def_Jpsi}.
Assuming further that the RH holds, we have
\[
  J_{k,\psi}(X,h)
  \ll
  h X^{1/k}  \log^2 \Bigl(\frac{2X}{h}\Bigr)
\]
uniformly for $ X^{1-1/k} \leq h  \leq X$.
\end{Lemma}

Theorem \ref{k-root-Saffari-Vaughan} is an immediate consequence of
Lemma~\ref{Lemma_da_h_a_deltax} and \eqref{fromThetatoPsiFinale}.
In its turn, Lemma \ref{Lemma_da_h_a_deltax} is a consequence of the
following result.

\begin{Lemma} 
\label{Lemma_deltax}
Let $k >0$ be a a real number and 
$\eps$ be an arbitrarily small positive constant.  
Assuming that \eqref{I-H} holds,
there exists a positive constant $c_1 = c_{1}(\eps)$, which does not
depend on $k$,  such that 
\begin{align}
\label{def-Jtilde}
\widetilde{J}_{k,\psi}(X,\delta)
&:=
\int_{X}^{2X}
\left(
\psi((x+\delta x)^{1/k})-\psi(x^{1/k}) -
((x+\delta x)^{1/k}-x^{1/k})
\right)^2
\dx x 
\\
\notag
&\ll
\delta^{2}X^{2/k+1}
\exp
\Big(
- c_{1}
\Big(
\frac{\log X}{\log \log X}
\Big)^{1/3}
\Big)
\end{align}
uniformly for $ X^{-2/(Ck)+\eps} \leq \delta \leq 1$. 
Assuming further that the RH holds, we have
\[
  \widetilde{J}_{k,\psi}(X,\delta)
  \ll
  \delta X^{1/k+1}  \log^2 \Bigl(\frac{2}{\delta}\Bigr)  
\]
uniformly for $ X^{-1/k} \leq \delta \leq 1$.

The same estimates hold if we insert $\theta$ in place of $\psi$
in the previous quantities.
\end{Lemma}
 
\begin{Proof}[Proof of Lemma \ref{Lemma_deltax}] 
We set $\Delta = (1 + \delta)^{1 / k} - 1$ in \eqref{def-Jtilde}
getting
\[
  \widetilde{J}_{k, \psi}(X, \delta)
  =
  \int_{X}^{2 X}
  \Bigl( \psi(x^{1/k} (1 + \Delta) ) - \psi( x^{1/k} ) - \Delta x^{1/k}
  \Bigr)^2 \, \dx x.
\]
Performing the substitution $y^k = x$ we get
\begin{align}
\notag
  \widetilde{J}_{k, \psi}(X, \delta)
  &=
  \int_{X^{1/k}}^{(2X)^{1/k}}
  \Bigl( \psi(y (1 + \Delta) ) - \psi(y) - \Delta y
  \Bigr)^2 \, k y^{k-1} \, \dx y \\
  &\asymp_k
\label{upp-bd-k<1}
  X^{1-1/k}
  J(X^{1/k}, \Delta).
\end{align}
In the unconditional case, by Lemma 5 of Saffari-Vaughan \cite{SaffariV1977}, 
for $\Delta > X^{-2 / (C k) + \eps}$ we have
\begin{align*}
  \widetilde{J}_{k,\psi}(X,\delta)
  &\asymp_{k}
  X^{1-1/k} \Delta^2 X^{3/k}
  \exp \Big( - c_{1}(\eps) \Big( \frac{\log X^{1/k}}{\log \log X^{1/k}}
                           \Big)^{1/3}
       \Big) \\
  &\asymp_{k}
  X^{2/k+1} \Delta^2 
  \exp \Big( - c_{2}(\eps, k) \Big( \frac{\log X}{\log \log X} \Big)^{1/3}
       \Big),
\end{align*}
where $c_{1}(\eps) > 0$ and we may take, essentially,
$c_{2}(\eps, k) = c_1(\eps) k^{-1 / 3}$.
For $k \in (0, 1)$ we are therefore allowed to take a value of $c_2$
which is independent of $k$.

In the conditional case, from \eqref{upp-bd-k<1} 
and by Lemma 5 of Saffari-Vaughan \cite{SaffariV1977}, we deduce
\[
  \widetilde{J}_{k,\psi}(X,\delta)
  \asymp_{k}
  X^{1-1/k} \Delta X^{2/k} \log\Bigl(\frac{2X^{1/k}}{\Delta}\Bigr)
  \asymp_{k}
  X^{1/k+1} \Delta \log\Bigl(\frac{2X}{\Delta}\Bigr).
\]
provided that $\Delta > X^{-1/k}$.
It is easy to see that
$\Delta = (1 + \delta)^{1/k} - 1 = (1/k) \delta + \Odip{k}{\delta^2}$,
so that $1 / \Delta = k / \delta + \Odip{k}{1} \ll_{k} 1 / \delta$.
Hence
\[
  \widetilde{J}_{k,\psi}(X,\delta)
  \asymp_{k}
  \begin{cases}
    X^{2/k+1} \delta^2 
    \exp \Big( - c_{2}(\eps,k) \Big( \dfrac{\log X}{\log \log X} \Big)^{1/3}
         \Big)
    &\text{unconditionally,} \\
    X^{1/k+1} \delta \log\Bigl(\dfrac{2X}{\delta}\Bigr)
    &\text{assuming RH.}
  \end{cases}
\]
A similar computation allows us to express the above bounds for
$\Delta$ in terms of $\delta$.
Skipping details for brevity, we may conclude that
Lemma~\ref{Lemma_deltax} holds true for $k \in (0, 1)$ (with a
constant $c_1 > 0$ that depends only on $\eps$) 
provided that
$\delta > X^{-2/(Ck)+2\eps}$ unconditionally, and that
$\delta > (1+\eps)kX^{-1/k}$ if we assume the RH. 

In the remaining range $k>1$ the previous proof gives a constant
$c_1$ which depends on $k$.
In fact, one can obtain Lemma~\ref{Lemma_deltax} in the full range for
$k$ and $c_1$ independent from $k$ following the proof of
Saffari-Vaughan \cite{SaffariV1977} (as in \cite{LanguascoS2012} for
the case $k=2$) but, since in the applications $k$ is usually bounded,
we omit such a proof here.

We finally remark that the estimates with $\theta$ in place of $\psi$ follow arguing
as in \eqref{fromThetatoPsiFinale} and using
the second part of  Lemma \ref{psi-theta-passaggio}.
\end{Proof}

\begin{Proof}
[Proof of Lemma \ref{Lemma_da_h_a_deltax}]
We follow the argument of \S 6 in Saffari-Vaughan \cite{SaffariV1977}. 
Let now  $2h \leq v \leq 3h$.
and define
\begin{equation} 
\label{diff-def}
\D_{k,\psi} (a,b)
=
\psi(a^{1/k})-\psi(b^{1/k}) - (a^{1/k}-b^{1/k}).
\end{equation}

To estimate $J_{k,\psi}(X,h)$ (defined in \eqref{def_Jpsi}), we first remark,
by \eqref{diff-def},
that
\begin{align}
h
J_{k,\psi} (X,h) 
& \ll
\int_{X}^{2X} 
\int_{2h}^{3h}  
\D^{2}_{k,\psi} (x+v,x)\
\dx v \
\dx x 
 +
 \int_{X}^{2X} 
\int_{2h}^{3h}  
\D^{2}_{k,\psi} (x+v,x+h)\
\dx v \
\dx x .
\label{int_da_studiare}
\end{align}
Setting $z=v-h, y = x+h$ and changing variables in the last integration,
the right-hand side of \eqref{int_da_studiare} becomes
\begin{align*}
 & \ll
\int_{X}^{2X} 
\int_{2h}^{3h} 
\D^{2}_{k,\psi} (x+v,x)\
\dx v \
\dx x  
 +
 \int_{X+h}^{2X+h} 
 \int_{h}^{2h} 
\D^{2}_{k,\psi} (y+z,y)\
\dx z \
\dx y.
\end{align*}
Since both the integrand functions are non-negative, we can extend the
integration ranges to get  
\begin{align*}
h
J_{k,\psi} (X,h)  
& \ll
\int_{X}^{2X+h} 
\int_{h}^{3h} 
\D^{2}_{k,\psi} (x+v,x)\
\dx v \
\dx x 
=
\int_{X}^{2X+h} 
x
\int_{h/x}^{3h/x} 
\D^{2}_{k,\psi} (x+\delta x,x)\
\dx \delta \ 
\dx x,
\end{align*}
where in the last step we made the change of variable $\delta = v/x$,
thus getting $\delta \geq h/x \geq X^{-2/(Ck) + \eps} $ as in 
the hypothesis of Lemma \ref{Lemma_deltax}.
Interchanging the integration order we obtain
\begin{align*}
h
J_{k,\psi}& (X,h)  
\ll
(X+h) 
\int_{h/(2X+h)}^{3h/X} 
\int_{X}^{2X+h} 
\D^{2}_{k,\psi} (x+\delta x,x)\
\dx x  \
\dx \delta.
\end{align*}
Finally, in the first case, i.e. assuming \eqref{I-H}, we use Lemma
\ref{Lemma_deltax} to get
\[
  J_{k,\psi}(X,h)
  \ll 
  h^{2} X^{2/k-1}
  \exp \Big( - c_{1} \Big( \frac{\log X}{\log \log X} \Big)^{1/3} \Big).
\]
Assuming RH, Lemma \ref{Lemma_deltax} implies
\begin{equation*}
J_{k,\psi}(X,h)
\ll
\frac{X+h}{h}
\int_{h/(2X+h)}^{3h/ X} 
\delta X^{1/k+1}  \log^2    \Bigl(\frac{2}{\delta}\Bigr)   \,
\dx \delta
\ll
h X^{1/k}  \log^2   \Bigl(\frac{2X}{h}\Bigr)  .
\end{equation*}
This concludes the proof of Lemma \ref{Lemma_da_h_a_deltax}.
\end{Proof}

The following elementary lemma is useful in passing 
from the $\theta$ to the $\psi$ function in Theorem \ref{k-root-Saffari-Vaughan}.
\begin{Lemma}
\label{psi-theta-passaggio}
Let $k > 0$ be a real number.
For $X^{1 - 1/k} \le h \le X$, we have
\begin{equation*}
  \int_X^{2 X}
    \Bigl(
    \psi\bigl( (x + h)^{1 / k} \bigr)
  -
  \psi\bigl( x^{1 / k} \bigr) 
  -
  \theta\bigl( (x + h)^{1 / k} \bigr)
  +
  \theta\bigl( x^{1 / k} \bigr) 
    \Bigr)^2 \, \dx x \\
  \ll
  h X^{1/k}.
\end{equation*}
Moreover, for  $X^{- 1/k} \le \delta \le 1$, we have
\begin{equation*}
 \int_X^{2 X}
    \Bigl(
    \psi\bigl( (x + \delta x)^{1 / k} \bigr)
  -
  \psi\bigl( x^{1 / k} \bigr) 
  -
  \theta\bigl( (x +  \delta x)^{1 / k} \bigr)
  +
  \theta\bigl( x^{1 / k} \bigr) 
    \Bigr)^2 \, \dx x \\
  \ll
   \delta  X^{1/k+1}.
\end{equation*}
\end{Lemma}

\begin{Proof}
Since $\psi(u) = \sum_{m=1}^{\log_2 u} \theta(u^{1/m})$, we have
\begin{align}
\notag
  &\psi\bigl( (x + h)^{1 / k} \bigr)
  -
  \psi\bigl( x^{1 / k} \bigr) 
  -
  \theta\bigl( (x + h)^{1 / k} \bigr)
  +
  \theta\bigl( x^{1 / k} \bigr) \\
  \label{decomp}
  &=
  \sum_{m=2}^{\log_2 (x^{1/k})}
    \Bigl( \theta\bigl( (x + h)^{1/mk} \bigr) - \theta\bigl( x^{1/mk} \bigr)
    \Bigr)
  +
  \sum_{m=\log_2 (x^{1/k})}^{\log_2 ((x+h)^{1/k})} \theta\bigl( (x + h)^{1/mk} \bigr).
\end{align}
Clearly, the last sum has at most $1 + 1 / k$ summands, which are
uniformly bounded.
 
Assume now that $h \in \bigl[X^{1 - 1 / k}, X^{1 - 1 / 2 k} \bigr]$ and denote
as $\Delta_k(X, h)$ the left-hand side of the inequality in the statement.  
Using \eqref{decomp}, we find that
\begin{align*}
  \Delta_k(X, h)
  &=
  \int_X^{2 X}
    \Bigl(
      \sum_{m = 2}^{\log_2(x) / k}
        \bigl( \theta \bigl( (x + h)^{1 / m k} \bigr)
               -
               \theta \bigl( x^{1 / m k} \bigr)
        \bigr)
      +
      \Odi{1}
    \Bigr)^2 \, \dx x \\
  &\ll
  \int_X^{2 X}
    \bigl( \theta \bigl( (x + h)^{1 / 2 k} \bigr)
           -
           \theta \bigl( x^{1 / 2 k} \bigr)
    \bigr)^2 \, \dx x \\
  &\qquad+
  \int_X^{2 X}
    \Bigl(
      \sum_{m = 3}^{\log_2(x) / k}
        \bigl( \theta \bigl( (x + h)^{1 / m k} \bigr)
               -
               \theta \bigl( x^{1 / m k} \bigr)
        \bigr)
    \Bigr)^2 \, \dx x
  +
  \Odi{X}.
\end{align*}
We only deal with the first term, the other one being similar and, in
fact, smaller.
We exploit the fact that the integrand is usually $0$ and small when
positive.
Let $M = \bigl\lceil X^{1 / 2 k} \bigr\rceil$ and
$N = \bigl\lfloor (2 X)^{1 / 2 k} \bigr\rfloor$.
Let $n \in [M, N]$; if $x \in [X, 2 X]$ satisfies both
$x^{1 / 2 k} < n$ and $(x + h)^{1 / 2k} \ge n$, then the integrand
$\bigl( \theta((x + h)^{1/2k}) - \theta(x^{1/2k}) \bigr)^2$ is
$\log^2 n$ if $n$ is a prime number, and vanishes otherwise.
The two inequalities above imply $n^{2 k} - h \le x < n^{2 k}$.
Summing up, for every prime $p \in [M, N]$ there is an interval of
length $h$ of values for $x$ such that the integrand does not vanish.
Hence
\begin{equation}
\label{first-range}
\Delta_k(X, h)
\ll
  h \sum_{p \in [M, N]} (\log p)^2
  \ll
  h X^{1 / 2 k}   \log X . 
\end{equation}

Let us consider now $h \in \bigl[X^{1 - 1 / 2k}, X \bigr]$.
For the terms with $m \ge 3$ in \eqref{decomp} we simply notice that
\begin{align*}
  \theta\bigl( (x + h)^{1/mk} \bigr) - \theta\bigl( x^{1/mk} \bigr)
  &\le
  \log\bigl( (x + h)^{1/mk} \bigr)
  \sum_{x^{1/mk} < n \le (x + h)^{1/mk}} 1 \\
  &\le
  \frac1{m k}
  \log(x + h)
  \Bigl( (x + h)^{1/mk} - x^{1/mk} + 1 \Bigr) \\
  &\le
  \frac1{m k}
  \log(x + h)
  \Bigl( \frac1{m k} h x^{1/mk - 1} + 1 \Bigr)
\end{align*}
by the mean-value theorem.
The number of such terms is at most $\log x$ and hence the total
contribution is bounded by an absolute constant times
$
  \log x \* \log \log x
  \bigl( h x^{1 / 3 k - 1} + 1 \bigr)$.
We now deal with the term
$\theta\bigl( (x + h)^{1/2k} \bigr) - \theta\bigl( x^{1/2k} \bigr)$.
We have
\begin{align*}
  \theta\bigl( (x + h)^{1/2k} \bigr) - \theta\bigl( x^{1/2k} \bigr)
  &\le
  \log\bigl( (x + h)^{1/2k} \bigr)
  \Bigl( \pi\bigl( (x + h)^{1/2k} \bigr) - \pi\bigl( x^{1/2k} \bigr)
  \Bigr) \\
  &\ll_k
  \log x
  \frac{h x^{1 / 2 k - 1} + 1}{\log(h x^{1 / 2 k - 1} + 1)}
\end{align*}
by the mean-value theorem again and the Brun-Titchmarsh inequality,
which we can use for $h \gg_k X^{1 - 1/2k}$. 
Squaring out and integrating the previous estimates we get, 
for every fixed $\eps>0$, that
\begin{equation}
\label{second-range}
  \Delta_k(X, h)
  \ll_{k,\eps}
  \begin{cases}
    h^{2} X^{1 / k - 1}        &\text{if $X^{1 - 1 / (2 k) + \eps} \le h \le X$,} \\
    h^{2} X^{1 /  k - 1} \log^{2} X &\text{if $X^{1 - 1 / (2 k)} \le h \le X^{1 - 1 / (2 k) + \eps}$.} 
  \end{cases}
\end{equation}  
The first part of Lemma \ref{psi-theta-passaggio} now follows from \eqref{first-range}-\eqref{second-range} by trivial computations.

The second part of Lemma \ref{psi-theta-passaggio} can be obtained
with a similar argument.
\end{Proof}

\renewcommand{\bibliofont}{\normalsize}

 \providecommand{\bysame}{\leavevmode\hbox to3em{\hrulefill}\thinspace}
\providecommand{\MR}{\relax\ifhmode\unskip\space\fi MR }
\providecommand{\MRhref}[2]{%
  \href{http://www.ams.org/mathscinet-getitem?mr=#1}{#2}
}
\providecommand{\href}[2]{#2}

\bigskip \noindent
\author{Alessandro LANGUASCO\\
Universit\`a di Padova\\
Dipartimento di Matematica\\
Via Trieste 63\\
35121 Padova, Italy\\
E-mail: languasco@math.unipd.it}

\bigskip \noindent
\author{Alessandro ZACCAGNINI \\
Universit\`a di Parma \\
Dipartimento di Matematica \\
Parco Area delle Scienze, 53/a \\ 
43124 Parma, Italy \\
E-mail: alessandro.zaccagnini@unipr.it}

\begin{thebibliography}{10}

\bibitem{BakerH1982}
R.~C. Baker and G.~Harman, \emph{Diophantine approximation by prime numbers},
  J. London Math. Soc. \textbf{25} (1982), 201--215.

\bibitem{BrudernCP1997}
J.~Br{\"u}dern, R.~J. Cook, and A.~Perelli, \emph{The values of binary linear
  forms at prime arguments}, {Proc. of Sieve Methods, Exponential sums and
  their Application in Number Theory} ({G. R. H. Greaves {\it et al}}, ed.),
  Cambridge U.P., 1997, pp.~87--100.

\bibitem{Gallagher1970}
P.~X. Gallagher, \emph{A large sieve density estimate near $\sigma=1$}, Invent.
  Math. \textbf{11} (1970), 329--339.

\bibitem{Ghosh1981}
A.~Ghosh, \emph{The distribution of $\alpha p^2$ modulo $1$}, Proc. London
  Math. Soc. \textbf{42} (1981), no.~2, 252--269.

\bibitem{Harman1991}
G.~Harman, \emph{Diophantine approximation by prime numbers}, J. London Math.
  Soc. \textbf{44} (1991), 218--226.

\bibitem{Huxley1972}
M.~N. Huxley, \emph{On the difference between consecutive primes}, Invent.
  Math. \textbf{15} (1972), 155--164.

\bibitem{LanguascoS2012}
A.~Languasco and V.~Settimi, \emph{On a {D}iophantine problem with one prime,
  two squares of primes and $s$ powers of two}, Acta Arithmetica \textbf{154}
  (2012), 385--412.

\bibitem{LanguascoZ2012c}
A.~Languasco and A.~Zaccagnini, \emph{A {D}iophantine problem with a prime and
  three squares of primes}, to appear in J. Number Theory,
  \url{http://arxiv.org/pdf/1206.0246}, 2012.

\bibitem{LanguascoZ2012e}
A.~Languasco and A.~Zaccagnini, \emph{A {Diophantine} problem with prime
  variables}, \url{http://arxiv.org/pdf/1206.0252}, Submitted, 2012.

\bibitem{Rieger1968}
G.~J. Rieger, \emph{\"{U}ber die {S}umme aus einem {Q}uadrat und einem
  {P}rimzahlquadrat}, J. Reine Angew. Math. \textbf{231} (1968), 89--100.

\bibitem{SaffariV1977}
B.~Saffari and R.~C. Vaughan, \emph{On the fractional parts of $x/n$ and
  related sequences. {II}}, Ann. Inst. Fourier \textbf{27} (1977), 1--30.

\bibitem{Tolev1992}
D.~Tolev, \emph{On a {D}iophantine inequality involving prime numbers}, Acta
  Arith. \textbf{61} (1992), 289--306.

\bibitem{Vaughan1974a}
R.~C. Vaughan, \emph{Diophantine approximation by prime numbers. {I}}, Proc.
  London Math. Soc. \textbf{28} (1974), 373--384.

\bibitem{Vaughan1974b}
R.~C. Vaughan, \emph{Diophantine approximation by prime numbers. {II}}, Proc.
  London Math. Soc. \textbf{28} (1974), 385--401.

\bibitem{Vaughan1997}
R.~C. Vaughan, \emph{The {Hardy}-{Littlewood} method}, second ed., Cambridge U.
  P., 1997.

\bibitem{Zaccagnini1998}
A.~Zaccagnini, \emph{Primes in almost all short intervals}, Acta Arith.
  \textbf{84.3} (1998), 225--244.

\end{thebibliography}
\end{document}